\documentclass[12pt]{article}

\usepackage[latin1]{inputenc}
\usepackage{amsmath,amssymb}
\usepackage{latexsym}

\newtheorem{theorem}{Theorem}[section]
\newtheorem{corollary}[theorem]{Corollary}

\newtheorem{proposition}[theorem]{Proposition}
\newtheorem{definition}[theorem]{Definition}
\newtheorem{assumption}[theorem]{Assumption}
\newtheorem{remark}[theorem]{Remark}
\numberwithin{equation}{section}

\def\proof{{\medskip\noindent {\bf Proof. }}}

\def\qed{{\hfill $\square$ \bigskip}}

\def\square{{\vcenter{\vbox{\hrule height.3pt
        \hbox{\vrule width.3pt height5pt \kern5pt
           \vrule width.3pt}
        \hrule height.3pt}}}}

  \def\sF {{\cal F}}
\def\sG {{\cal G}} \def\sH {{\cal H}} 
  \def\sL {{\cal L}}
\def\sM {{\cal M}}  
  \def\sR {{\cal R}}
\def\sS {{\cal S}}  
\def\sV {{\cal V}}

\def\wt{\widetilde}
\def\wh{\widehat}
\def\ol{\overline}
\def\E{{\mathbb E}}
\def\P{{\mathbb P}}
\def\norm#1{{\Vert #1 \Vert}}
\def\del{{\partial}}
\def\lam{{\lambda}}

\def\bee{\begin{equation}}
\def\eee{\end{equation}}

\def\R{{\mathbb R}}
\def\E{{{\mathbb E}\,}}
\def\P{{\mathbb P}}

\def\F{{\cal F}}

\def\lam{{\lambda}}

\def\al{{\alpha}}
\def\grad{{\nabla}}
\def\proof{{\medskip\noindent {\bf Proof. }}}

\def\qed{{\hfill $\square$ \bigskip}}

\def\eps{\varepsilon}

\def\norm#1{\Vert #1 \Vert}

 \def\qq {\qquad}
\def\del{{\partial}}
\def\wt{\widetilde}
\def\ol{\overline}

\def\wh{\widehat}

\def\ni{\noindent }
\def\ms{\medskip}

\def\real{{\rm Re\,}}

\def\square{{\vcenter{\vbox{\hrule height.3pt
        \hbox{\vrule width.3pt height5pt \kern5pt
           \vrule width.3pt}
        \hrule height.3pt}}}}

\def\bdot{{\textstyle{.}}}

\def\tlint{{- \kern-0.85em \int \kern-0.2em}}  % for textstyle
\def\dlint{{- \kern-1.05em \int \kern-0.4em}}  % for displays

  \def\sF {{\cal F}}
\def\sG {{\cal G}} \def\sH {{\cal H}} 
  \def\sL {{\cal L}}
\def\sM {{\cal M}}  
  \def\sR {{\cal R}}
\def\sS {{\cal S}}  
\def\sV {{\cal V}}

\def\nn{{\nonumber}}

\begin{document}

\title{The martingale problem for a class of stable-like processes}

\author{Richard F. Bass\footnote{Research partially supported by NSF grant
DMS-0601783.}  \; and Huili Tang}

%\date{January 5,  2005}

\maketitle

\begin{abstract}  
\noindent Let $\al\in (0,2)$ and consider the operator
$$\sL f(x) =\int [f(x+h)-f(x)-1_{(|h|\leq 1)} \grad f(x)\cdot h]
\frac{A(x,h)}{|h|^{d+\al}}\, dh, $$
where the $\grad f(x)\cdot h$ term is omitted if $\al<1$. We consider the
martingale problem corresponding to the operator $\sL$ and 
under mild conditions on the function $A$ prove that there exists a unique solution.


\end{abstract}

\section{Introduction}\label{sect:intro}

A stable-like process is a pure jump process where the jump
intensity kernel is comparable in some sense to that  of one or more stable processes. The 
term was introduced in \cite{Discl}  for processes whose
associated operators were of the form
$$\int_\R [f(x+h)-f(x)-1_{(|h|\leq 1)} \grad f(x)\cdot h]
\frac{dh}{|h|^{1+\al(x)}},$$ 
 and the use of the term
was extended in \cite{CK} to refer to symmetric Markov processes whose
jump kernels $J(x,y)$ were comparable to $|x-y|^{-d-\al}$ for a fixed $\al$.

In this paper we fix $\al\in (0,2)$.  For $\al\in [1,2)$ we consider jump processes
associated to the operator
\bee\label{Op1}\sL f(x)=\int [f(x+h)-f(x)-1_{(|h|\leq 1)} \grad f(x)\cdot h]
\frac{A(x,h)}{|h|^{d+\al}}\, dh, 
\eee 
and for $\al\in (0,1)$ associated to the operator
\bee\label{Op2}\sL f(x)=\int [f(x+h)-f(x)]
\frac{A(x,h)}{|h|^{d+\al}}\, dh, 
\eee 
where $A(x,h)$ is bounded above and below by positive constants not
depending on $x$ or $h$. 
For the domain of $\sL$ we take the class of $C^2$ functions such that the 
function and its first and second partial derivatives are bounded.
These jump processes, when at a point $x$, jump
to $x+h$ with intensity given by $A(x,h)|h|^{-d-\al}$. These processes
stand in the same relationship to symmetric stable processes of
index $\al$ as uniformly elliptic operators in non-divergence form
do to Brownian motion. 

For $\al\geq 1$
the $\grad f(x)\cdot h$ term is needed to guarantee convergence of the integral, 
while for $\al<1$ the $\grad f(x)\cdot h$ term cannot be present, or else the
jumps of the process will be dominated by the drift. 

Processes corresponding to $\sL$ given by \eqref{Op1} or \eqref{Op2}
were considered in \cite{Jpharn} and \cite{SV}, where Harnack inequalities
and regularity of harmonic functions were proved. It is natural to
ask whether there exists a process corresponding to $\sL$, and if so,
is there only one. 

We view this question as a martingale
problem. Let $\Omega=D([0,\infty))$, the set of paths that are right
continuous with left limits, endowed with the Skorokhod topology.
Set $X_t(\omega)=\omega(t)$ for $\omega\in \Omega$, define $\theta_t:
\Omega\to \Omega$ by $\theta_t(\omega)(s)=\omega(s+t)$, and let
$\sF_t$ be the right continuous filtration generated by the process $X$.
A probability measure $\P$ is a solution
to the  martingale problem for $\sL$ started at $x$ if $\P(X_0=x)=1$ and
$f(X_t)-f(X_0)-\int_0^t \sL f(X_s)\, ds$ is a martingale whenever 
$f$ is a $C^2$ function such that $f$ and its first and second partial
derivatives are bounded. The question to be answered is the existence
and uniqueness of a solution to the martingale problem for $\sL$.

Our results on existence are merely an application of techniques used
in \cite{Upjump} and the novelty in the current paper is a sufficient condition
for uniqueness. Let $\eta>0$ and set
$$\psi_\eta(x)=(1+\log^+(1/x))^{1+\eta}, \qq x>0.$$
We require continuity in $x$ of the function $A(x,h)$, with 
more continuity the smaller $h$ is. More specifically, let
$$\ol A(x,h)=A(x,h)\psi_\eta(|h|).$$

Our main assumption is that $\ol A(x,h)$ be continuous in $x$, uniformly in
$h$. We assume

\begin{assumption}\label{A.0}
(a) There exist $c_1, c_2>0$ such that $c_1\leq A(x,h)\leq c_2$ for all $x$
and $h$.

(b)  There exists $\eta>0$  such that for every $y\in \R^d$ and every $b>0$
$$\lim_{x\to y}\sup_{|h|\leq b} |\ol A(x,h)-\ol A(y,h)|=0.$$
\end{assumption}

\ni Part (a) of Assumption \ref{A.0} may be regarded as the jump process equivalent
of uniform ellipticity.

We then have 

\begin{theorem}\label{main}
Suppose Assumption \ref{A.0} holds and $x_0\in \R^d$. Then there is one and only one
solution to the martingale problem for $\sL$ started at 
$x_0$.
\end{theorem}

As we alluded to above, existence is already known and can in fact be
proved under slightly weaker hypotheses. 
Some other generalizations are possible;
see Remarks \ref{rem:psi} and \ref{rem:largejumps}.
Our theorem also extends some of the results obtained
in  \cite{komatsu};
see Remark \ref{rem:kom}. 

We do not know if our theorem is still true if $\ol A$ is replaced by
$A$ in Assumption \ref{A.0}. We point out that uniqueness for
the martingale problem for jump 
processes  does not always hold; see \cite[Section 6]{Nlharn}.

In the next section we establish some estimates. An approximation is given
in Section \ref{sect:approx} and Theorem \ref{main}
is proved in Section \ref{sect:exist}.

\section{Estimates}\label{sect:est}

Let $B(x,r)=\{y\in \R^d: |y-x|<r\}$. 
Let $C^k$ be the functions which are $k$ times continuously differentiable,
$C^k_b$ the elements of $C^k$ such that the function and its partial derivatives
up to order $k$ are bounded, and $C^k_K$ the functions in $C^k$ that have compact support.
We use the probabilist's version of the Fourier transform:
$$\wh f(u)=\int e^{iu\cdot x} f(x)\, dx.$$
For processes whose paths are right continuous  with left limits,
we set $X_{t-}=\lim_{s<t, s\to t} X_s$ and $\Delta X_t=X_t-X_{t-}$.
We use the letter $c$ with or without subscripts to denote
constants whose value is unimportant and may change from line to line.

We suppose throughout the remainder of the paper that Assumption \ref{A.0}
holds.

\begin{definition}\label{D.SMF}
{\rm We say a collection $\{\P^x\}$ of probability measures is a strong Markov family of solutions
to the martingale problem for $\sL$ if for each $x\in \R^d$, $\P^x$ is a solution
to the martingale problem for $\sL$ started at $x$ and in addition
the strong Markov property holds: for any finite stopping time $T$, 
any $Y$ bounded and $\F_\infty$-measurable, and any $x\in \R^d$,
$$\E^x[Y\circ \theta_T\mid \F_T]=\E^{X_T}[Y], \qq \P^x-\mbox{a.s.}$$
}
\end{definition}

\begin{proposition}\label{P3.1A}
Suppose $r<1$, $x \in \R^d$, $\tau_r=\inf\{t: |X_t-x|\geq r\}$, and $\P$
is a solution to the martingale problem for
$\sL$ started at $x$. There exists $c_1$ not depending on $x$ such that
$$\P(\sup_{s\leq t} |X_s-x|\geq r)\leq c_1t/r^2, \qq t>0.$$
\end{proposition}

\proof Let $f:\R^d\to [0,1]$ be a  $C^2$  function such that $f(0)=0$ and
$f(y)=1$ if $|y|>1$. Let $f_{rx}(y)=f((y-x)/r)$. There exists
a constant $c$ such that  the first derivatives of $f_{rx}$ are bounded
by $c/r$ and the second derivatives are bounded by $c/r^2$.
By Taylor's theorem,
$$|f_{rx}(z+h)-f_{rx}(z)-\grad f_{rx}(z)\cdot h|\leq c|h|^2/r^2$$
and
$$|f_{rx}(z+h)-f_{rx}(z)|\leq c |h|/r.$$

Suppose $\al\geq 1$. Then
\begin{align*}
|\sL f_{rx}(z)|&\leq \int_{|h|\leq r} |f_{rx}(z+h)-f_{rx}(z)-\grad
f_{rx}(z)\cdot h|
\frac{A(x,h)}{|h|^{d+\al}}\, dh\\
&\qq +\int_{1\geq |h|>r} |f_{rx}(z+h)-f_{rx}(z)-\grad f_{rx}(z)\cdot h|
\frac{A(x,h)}{|h|^{d+\al}}\, dh\\
&\qq +\int_{|h|>1} |f_{rx}(z+h)-f_{rx}(z)|
\frac{A(x,h)}{|h|^{d+\al}}\, dh\\
&\leq \frac{c}{r^2} \int_{|h|\leq r}\frac{|h|^2}{|h|^{d+\al}}\, dh
+\frac{c}{r} \int_{|h|>r} \frac{|h|}{|h|^{d+\al}}\, dh\\
&\leq c r^{-\al}.
\end{align*}
Therefore by Doob's optional stopping theorem
\begin{align*}
\P(\tau_r\leq t)&\leq \E f_{rx}(X_{\tau_r\land t})-f_{rx}(x)\\
&=\E \int_0^{\tau_r\land t} \sL f_{rx}(X_s)\, ds\\
&\leq c t/r^\al.
\end{align*}

The case $\al<1$ is similar.
\qed

\begin{proposition}\label{P3.3} If $f\in C^2_b$, then $\sL f$ is continuous.
\end{proposition}

\proof Let $\eps>0$ and suppose that $\al\geq 1$, the case when
$\al<1$ being very similar. Let $\delta\in (0,1)$ and write
\begin{align*}
\sL f(x)&=\int_{|h|\leq \delta} [f(x+h)-f(x)-\grad f(x)\cdot h]
\frac{A(x,h)}{|h|^{d+\al}}\, dh\\
&\qq +\int_{\delta< |h|\leq 1} [f(x+h)-f(x)-\grad f(x)\cdot h]
\frac{A(x,h)}{|h|^{d+\al}}\, dh\\
&\qq +\int_{1< |h|\leq \delta^{-1}} [f(x+h)-f(x)]
\frac{A(x,h)}{|h|^{d+\al}}\, dh\\
&\qq +\int_{\delta^{-1}< |h|} [f(x+h)-f(x)]
\frac{A(x,h)}{|h|^{d+\al}}\, dh.
\end{align*}
The first term is bounded by
$$c\int_{|h|\leq \delta} \frac{|h|^2}{|h|^{d+\al}}\, dh,$$
where $c$ depends on $f$. This
is less than $\eps$ if $\delta$ is sufficiently small. 
The fourth term is bounded by
$$c\int_{|h|> \delta^{-1}} \frac{dh}{|h|^{d+\al}},$$
where again $c$ depends on $f$. This
will also be less than $\eps$ if $\delta$ is sufficiently small. 
The
second and third terms  are continuous in $x$ by dominated convergence and the
continuity of $A(x,h)$ in $x$.
\qed

\begin{proposition}\label{P3.1B}
Suppose  $\{\P^x\}$ is a strong Markov family of
solutions to the martingale problem for $\sL$.
Let $x_0\in \R^d$, suppose $r<1$, and $\tau_r=\inf\{t: |X_t-x_0|>r\}$.
\begin{description}
\item{(a)} If $\eps>0$, there exists $c_1$ (depending on $\eps$)
such that
$$\inf_{z\in B(x_0,(1-\eps)r)} \E^z \tau_r\geq c_1r^{-\al}.$$
\item{(b)} There exists $c_2$ such that
$$\sup_{z} \E^z \tau_r\leq c_2 r^\al.$$
\end{description}
\end{proposition}

\proof
The proof consists of  minor modifications to the proofs of
\cite[Lemmas 3.2 and 3.3]{Upjump}.
\qed

\begin{proposition}\label{P3.21}
Let $\P$ be a solution to the martingale problem for $\sL$
started at some point $x_0$.
If $B$ and $C$ are Borel sets whose closures are disjoint, then
$$\sum_{s\leq t} 1_B(X_{s-})1_C(X_s)-\int_0^t 1_B(X_s)\int_C \frac{A(X_s, u-X_s)}{|u-X_s|^{d+\al}}\, du\, ds$$
is a martingale with respect to $\P$.
\end{proposition}

\proof Suppose $B$ and $C$ are disjoint compact sets, $f\in C^2_b$ is 0 on $B$ and 1 on $C$, and $\grad f$ is 0 on $B$.
Then
$$f(X_t)-f(X_0)=M_t+\int_0^t \sL f(X_s)\, ds,$$
where $M_t$ is a martingale. It follows that  $\int_0^t 1_B(X_{s-})\, dM_s$ is
also a martingale. By Ito's formula
$$f(X_t)-f(X_0)=\int_0^t \grad f(X_{s-})\cdot dX_s +\sum_{s\leq t} 
[f(X_s)-f(X_{s-})-\grad f(X_{s-})\cdot \Delta X_s].$$
Hence
\begin{align}
\int_0^t 1_B(X_{s-}) &\grad f(X_{s-})\cdot \Delta X_s\label{E3.9}\\
&+\sum_{s\leq t} 1_B(X_{s-})[f(X_s)-f(X_{s-})-\grad f(X_{s-})\cdot \Delta X_s]
\nn\\
&-\int_0^t 1_B(X_{s-}) \sL f(X_s)\, ds\nn
\end{align}
is a martingale. Since $f\in C^2$ and both $f$ and $\grad f$ are 0 on $B$,
the first term  of \eqref{E3.9} is equal to 
0 and the second term of \eqref{E3.9}  is
$$\sum_{s\leq t} 1_B(X_{s-}) f(X_s).$$
We have
\begin{align*}
1_B(x)\sL f(x)&=1_B(x)\int [f(x+h)-f(x)-1_{(|h|\leq 1)}\grad f(x)\cdot h] \frac{A(x,h)}{|h|^{d+\al}}\, dh\\
&=1_B(x)\int f(x+h)\frac{A(x,h)}{|h|^{d+\al}}\, dh\\
&=1_B(x) \int f(u) \frac{A(x,u-x)}{|u-x|^{d+\al}}\, du.
\end{align*}
Putting this in (\ref{E3.9}), and using the fact that $X_s$ differs from $ X_{s-}$ on a set
of times having Lebesgue measure 0, the last term in \eqref{E3.9} is
$$\int_0^t 1_B(X_s)\int f(u)\frac{A(X_s,u-X_s)}{|u-X_s|^{d+\al}}\, du\, ds.$$
Our result follows by using a  limit argument.
\qed

\begin{proposition}\label{P3.4} Suppose  $\{\P^x\}$ is a strong
Markov family of solutions to the martingale
problem for $\sL$. 
Suppose $g$ is bounded and $\lam>0$. Let
$$S_\lam g(x)=\E^x \int_0^\infty e^{-\lam t} g(X_t)\, dt.$$
Then $S_\lam g$ is H\"older continuous in $x$.
\end{proposition}

\proof The proof follows by \cite[Theorem 4.3]{Jpharn} and the arguments
leading up to it. 
See also \cite{SV}.
\qed

Let 
$\sM^z$ be the operator on $C^2_b$ functions defined by
\bee\label{DL0}
\sM^z f(x)=\int [f(x+h)-f(x)-\grad f(x)\cdot h1_{(|h|\leq 1)}]
\frac{A(z,h)}{|h|^{d+\al}}\, dh,
\eee
where  the $\grad f(x)\cdot h$
term is missing if $\al<1$. Let $R^z_\lam$ be the resolvent
for the L\'evy process whose infinitesimal generator is  $\sM^z$
and let $P_t^z$ be the corresponding transition operator.

\begin{proposition}\label{P3.6}
With $\sM^z$ as above,
$$\real (\wh \sM^z (u))\leq 
    \begin{cases}
     0,& \qq |u|\leq 1;\\
     -c|u|^\al, &\qq |u|>1.
     \end{cases}
$$
\end{proposition}

\proof 
$$\wh \sM^z (u)=\int [e^{iu\cdot h}-1-iu\cdot h1_{(|h|\leq 1)}]
\frac{A(z,h)}{|h|^{d+\al}}\, dh,$$
with the $iu\cdot h1_{(|h|\leq 1)}$ term missing if $\al<1$. So
$$-\real(\wh \sM^z(u))=\int[1-\cos(u\cdot h)]
\frac{A(z,h)}{|h|^{d+\al}}\, dh,$$
and the assertion in the case $|u|\leq 1$ is immediate.

If $|u|>1$, setting $u=rv$, where $|v|=1$ and $r\in (1,\infty)$,
\begin{align*}
-\real(\wh \sM^z(u))&\geq c\int_{|h|\leq 1}
[1-\cos(u\cdot h)]\frac{1}{|h|^{d+\al}}\, dh\\
&=c\int_{|h|\leq 1}[1-\cos(r(v\cdot h))]\frac{1}{|h|^{d+\al}}\, dh\\
&=\frac{c}{r^\al}\int_{|h|\leq r}
[1-\cos(v\cdot h)]\frac{1}{|h|^{d+\al}}\, dh\\
&\geq\frac{c}{r^\al}\int_{|h|\leq 1}
[1-\cos(v\cdot h)]\frac{1}{|h|^{d+\al}}\, dh,
\end{align*}
using a change of variables. The integral in the last line is bounded
below by
a constant, using rotational invariance, and the result follows on noting
$r=|u|$.
\qed

\begin{corollary}\label{C3.7}
If $p^z(t,x,y)=\ol p^z_t(x-y)$ is the transition density for the L\'evy process with
generator $\sM^z$, then for each $t$, $\sup_z \norm{\ol p^z_t}_2\leq c(t)<\infty$.

Moreover, if $r^z_\lam(x)=\int_0^\infty e^{-\lam t} \ol p^z_t(x)\, dt$, then 
$$|\wh {r^z}_\lam (u)|\leq \frac{c}{\lam+|u|^\al}.$$
\end{corollary}

\proof The Fourier transform of $\ol p^z_t$ is $e^{t\wh \sM^z(u)}$,
and so $$|e^{t\wh \sM^z(u)}|=e^{t\real(\wh \sM^z(u))}.$$
With the estimates from Proposition \ref{P3.6},
this is less than or equal to 1 if $|u|\leq 1$ and less than or
equal to $e^{-ct|u|^\al}$ if $|u|>1$, where $c$ does not depend on 
$z$. So the Fourier transform of
$\ol p^z_t$ is in $L^2$, hence $\ol p^z_t$ is in $L^2$ by Plancherel's theorem,
with a bound not depending on $z$.

Now
$$|\wh {r^z}_\lam (u)|=\Big|\frac{1}{\lam-\wh \sM^z(u)}\Big|
\leq \frac{1}{\real(\lam-\wh \sM^z(u))}.$$
This is less than or equal to $1/\lam$ if $|u|\leq 1$ and less than
$1/(\lam+c|u|^\al)$ if $|u|>1$.
\qed

\begin{proposition}\label{P3.5} 
If $f\in L^2$,
$\norm{R^z_\lam f}_2\leq \frac1{\lam}\norm{f}_2.$
\end{proposition}

\proof By Corollary \ref{C3.7}, $P^z_t$ has a density $p^z(t,x,y)
=\ol p^z_t(x-y)$ for some function
$\ol p^z_t$ in $L^1$. Then
$$\int \ol p^z_t(x-y)\, dx=\int \ol p^z_t(x-y)\, dy=1$$
by a change of variables. Hence $\norm{P^z_tf}_2\leq \norm{f}_2$ by
\cite[Theorem IV.5.1]{PTA}. We now apply  Minkowski's inequality for integrals.
\qed

\begin{proposition}\label{P3.8}
Let $R^z_\lam$ be as above, $f\in L^2\cap C^2_K$. 
\begin{description}
\item{(a)} If $\al<1$ and $g(x)=R^z_\lam f(x+h)-R^z_\lam f(x)$, then
$$\norm{g}_2\leq c|h|^\al \norm{f}_2.$$
\item{(b)} If $\al\in (0,2)$, then
$$\norm{g}_2\leq \frac{c}{\lam} \norm{f}_2.$$
\item{(c)} If $\al\in [1,2)$ and 
$$G(x)=R^z_\lam f(x+h)-R^z_\lam f(x)-\grad R^z_\lam f(x)\cdot h,$$
then
$$\norm{G}_2\leq c|h|^\al \norm{f}_2.$$
\item{(d)} If $\al\in [1,2)$, then
$$\norm{G}_2\leq c\Big(\frac{1}{\lam}+|h|\Big)\norm{f}_2.$$
\end{description}
\end{proposition}

\proof
First of all, if $f\in L^2\cap C^2_K$, then $R^z_\lam f\in L^2\cap C_b^2$ by 
Proposition \ref{P3.5} and translation invariance.
So $\grad R^z_\lam f$ is well defined. By translation invariance,
$\frac{\del R^z_\lam f}{\del x_i}=R^z_\lam (\frac{\del f}{\del x_i})$, and
$\frac{\del f}{\del x_i}\in C^1_K\subset L^2$, so $R^z_\lam(\frac{\del f}{\del x_i})
\in C^1_b$, and is in $L^2$. Therefore to prove the proposition it suffices
to look at Fourier transforms and to use Plancherel's theorem.

(a) We have
$$\wh g(u)=\wh f(u)\wh {r^z}_\lam(u)[e^{iu\cdot h}-1],$$
so using Corollary \ref{C3.7} 
$$|\wh g(u)|\leq \frac{c|\wh f(u)|}{\lam+|u|^\al}|h|^\al |u|^\al
\leq c|\wh f(u)|\, |h|^\al.$$
Therefore
$$\norm{\wh g}_2\leq c|h|^\al \norm{\wh f}_2,$$
and the result follows by Plancherel's theorem.

(b) As in (a), but using $|e^{iu\cdot h}-1|\leq 2$, we have
$$|\wh g(u)|\leq \frac{2c|\wh f(u)|}{\lam},$$
and we use Plancherel's theorem as in (a).

(c) $$\wh G(u)=\wh f(u)\wh {r^z}_\lam(u)[e^{iu\cdot h}-1-iu\cdot h].$$
Now
\begin{align*}
|e^{iu\cdot h}-1-iu\cdot h|&=\Big|\int_0^{u\cdot h} [ie^{is}-i]\, ds\Big|\\
&\leq c\int_0^{|u\cdot h|} |s|^{\alpha-1}\, ds\\
&\leq c |u\cdot h|^{\al}.
\end{align*}
Hence
$$|\wh G(u)|\leq c\frac{|\wh f(u)|}{\lam +|u|^\al}|u|^\al|h|^\al
\leq c|h|^\al |\wh f(u)|.$$

(d)
Similarly to the proofs of (b) and (c),
$$|\wh G(u)|\leq c\frac{|\wh f(u)|}{\lam+|u|^\al}(2+|u\cdot h|).$$
If $|u|\leq 1$, then 
$$|\wh G(u)|\leq \frac{c}{\lam} |\wh f(u)|(2+|h|).$$
On the other hand, if $|u|>1$, then since $\al\geq 1$ and
$$\frac{|u\cdot h|}{\lam+|u|^\al}\leq \frac{|u|\, |h|}{|u|^\al}\leq |h|,$$
we have
$$|\wh G(u)|\leq c|h| \, |\wh f(u)|.$$
\qed

Using this proposition we can extend the definition of the functions
$g,G$ and extend the above estimates to every $f\in L^2$.

\section{Approximation}\label{sect:approx}

A key step in the uniqueness proof is to get a bound on the resolvent
for an arbitrary solution to the martingale problem for $\sL$.
We do that by an approximation procedure.

We begin with

\begin{definition}\label{defPPP}
{\rm Let $(S,\sS,\lam)$ be a measure space, where $\lam$ is a $\sigma$-finite
measure.
A random measure $\mu([0,t]\times A)(\omega)$ is a Poisson point process
with intensity measure $\lam$
if

(a) whenever $\lam(A)<\infty$, $N_t(A)=\mu([0,t]\times A)$ is a  Poisson
process with intensity $\lam(A)$ and

(b) If $n\geq 1$ and $A_1, \ldots, A_n$ are disjoint with 
$\lam(A_1), \ldots, \lam(A_n)<\infty$
for each $i$, then the processes $N_t(A_i)$, $i=1, \ldots, n$, are
independent.
}
\end{definition}

\begin{proposition}\label{PA1}
Suppose whenever $\lam(A)<\infty$, $N_t(A)$ is a process starting at
0 with paths that are right continuous and left limits that are constant
except for jumps that are of size one. If $N_t(A)-\lam(A)t$ is
a martingale for each such $A$ and $N_t(A)$ and $N_t(B)$ have no jumps
in common when $A$ and $B$ are disjoint, then 
$\mu([0,t]\times A)=N_t(A)$ is a Poisson point process.
\end{proposition}

\proof Property (a) of the definition of Poisson point process follows
from a very slight modification of \cite[III.T12]{Meyer76}. In addition, that theorem shows that 
$\sigma(N_t(A)-N_s(A): A \in \sS\})$
is independent of $\F_s$. 

We next prove that if $A_i$, $i=1, \ldots, n$,  are disjoint sets of finite $\lam$-measure
and $t_0>0$, then 
\bee\label{EA1}
N_{t_0}(A_1), \ldots, N_{t_0}(A_n) \mbox{ are independent random variables.}
\eee

To prove \eqref{EA1}, we do the case when $n=2$, the general case being
very similar. Let $u_1, u_2$ be two reals and define 
$$M_t^j=\exp\Big(iu_j N_{t\land t_0}(A_j)-\lam(A_j)(t\land t_0)(e^{iu_j}-1)\Big), \qq j=1,2.$$
Because $N_t(A_j)$ is a Poisson process with intensity $\lam(A_j)$,
each $M_t^j$ is a martingale with $M_0^j=0$.

Since $N_t(A_1)$ and  $N_t(A_2)$ have no jumps in common and are non-de\-crea\-sing, the quadratic
variation process $[N_\bdot(A_1), N_\bdot(A_2)]_t$ is zero. So by
Ito's product formula,
$$M_\infty^1M_\infty^2=M_0^1M_0^2 +\int_0^\infty M^1_{s-}\, dM^2_s
+\int_0^\infty M^2_{s-}\, dM^1_s,$$
or
$\E[M_\infty^1M_\infty^2]=1.$ It follows that
\begin{align*}
\E\Big[e^{iu_1N_{t_0}(A_1)} e^{iu_2N_{t_0}(A_2)}\Big]&=
e^{\lam(A_1)t_0(e^{iu_1}-1)}e^{\lam(A_2)t_0(e^{iu_2}-1)}\\
&=\E\Big[e^{iu_1N_{t_0}(A_1)}\Big]\E\Big[e^{iu_2N_{t_0}(A_2)}\Big].
\end{align*}
This holds for every $u_1, u_2$, so $N_{t_0}(A_1)$ and $N_{t_0}(A_2)$
are independent.

A very similar argument shows that if $0<s_0<t_0$, then
$N_{t_0}(A_1)-N_{s_0}(A_1), \ldots, N_{t_0}(A_n)-N_{s_0}(A_n)$ are
independent random variables. This and the independence of 
$\sigma(N_t(A)-N_s(A): A\in \sS)$ from $\F_s$ implies part (b)
of Definition \ref{defPPP}.
\qed

We next construct a function $F(x,u)$ such that for every Borel set $B$
and every $x\in \R^d$
\bee\label{Eequiv}
\int 1_B(F(x,u))\frac{du}{|u|^{d+\al}} =\int_B \frac{A(x,h)}{|h|^{d+\al}}\, dh.
\eee
Such constructions are known (see \cite{LM} or \cite{Jacod}), but
we want our $F$ to be continuous in $x$ as well, and the existing
constructions do not necessarily possess this property. (When $d>1$, 
$F$ satisfying \eqref{Eequiv} are by no means unique.) We will do the case $d=2$ for simplicity of
notation, but the idea for higher dimensions is essentially the same.

Fix $x$. We define $F$ for $u$ in the first quadrant, and the other quadrants
are done similarly. Set $r_0=\infty$ and choose $r_1>r_2>\cdots >0$ such that
$$\int_{[r_{i+1},r_i)\times [0,\infty)} \frac{A(x,h)}{|h|^{d+\al}}\, dh=2^{-1}, \qq i=0,1,\ldots$$
For each strip
$[r_{i+1},r_i)\times [0,\infty)$, let $s_0=\infty$ and choose
$s_1>s_2>\cdots > 0$ such that
$$\int_{[r_{i+1},r_i)\times [s_{j+1},s_j)} \frac{A(x,h)}{|h|^{d+\al}}\, dh=2^{-2},\qq j=0,1,\ldots$$
Let $\sR_1=\sR_1(x)$ be the collection of such rectangles. 
Note that each rectangle in $\sR_1$ has the same mass with respect to
the measure
$A(x,h)/|h|^{d+\al}\, dh$, but the rectangles are not congruent in shape.

Set $v_0=\infty$ and choose $v_1>v_2>\cdots >0$ such that
$$\int_{[v_{i+1},v_i)\times [0,\infty)} \frac{du}{|u|^{d+\al}}=2^{-1},\qq i=0,1,\ldots$$
For each strip $[v_{i+1},v_i)\times [0,\infty)$, let $w_0=\infty$ and choose
$w_1>w_2>\cdots > 0$ such that
$$\int_{[v_{i+1},v_i)\times [w_{j+1},w_j)} \frac{du}{|u|^{d+\al}}=2^{-2},\qq j=0,1,\ldots$$
Let $\sV_1=\sV_1(x)$ be the collection of such rectangles.
 
Let $\Gamma_1$ be the map from $\sV_1$ to $\sR_1$ taking the element
$[v_{i+1},v_i)\times [w_{j+1},w_j)$ of $\sV_1$ 
to the element
$[r_{i+1},r_i)\times [s_{j+1},s_j)$ of $\sR_1$. 

If $[r, r')\times [s,s')$ is an element of $\sR_1$, choose 
$r''$ such that
$$\int_{[r,r'')\times [s,s')} \frac{A(x,h)}{|h|^{d+\al}}\, dh=2^{-3},$$
and then $s'', s'''$ such that the integrals of $A(x,h)/|h|^{d+\al}$
over $[r,r'')\times [s,s'')$
and over $[r'',r')\times [s,s''')$ are both equal to $2^{-4}$. 
We put the 4 rectangles $[r,r'')\times [s,s'')$, $[r,r'')\times [s'',s')$,
$[r'',r')\times [s,s''')$, and $[r'',r')\times [s''', s')$ into $\sR_2=\sR_2(x)$
and do this for each rectangle in $\sR_1$. We divide each rectangle of 
$\sV_1$ similarly into 4 rectangles of mass $2^{-4}$ with respect to
the measure $du/|u|^{d+\al}$ and let $\sV_2$ be the collection of such
subrectangles. We define the map $\Gamma_2$ from $\sV_2$ into $\sR_2$
that takes a rectangle of $\sV_2$ into the corresponding rectangle of $\sR_2$. 

We continue by dividing each rectangle of $\sR_2$ and $\sV_2$ into 4 subrectangles, and so on.
Now define $F_m(x,u):(0,\infty)^2\to (0,\infty)^2$ by setting $F_m(x,u)$
to be the lower left hand point of $\Gamma_m(U)$ if $u\in U\in \sV_m$.
Recall $x$ is fixed.
It is easy to check that
$F_m(x,u)$ converges uniformly over $u$ in compact subsets of $(0,\infty)^2$,
and if we call the limit $F(x,u)$, then $u\to F(x,u)$ will be one-to-one.
The construction shows that the equality
\eqref{Eequiv} holds if $m\geq 1$ and $B\in \sV_m$, and therefore
it holds for every Borel set contained in $(,\infty)^2$.

Moreover, if $x_0\in \R^d$ is fixed and $m\geq 1$, then by taking $x$
sufficiently close to $x_0$ the boundaries of each rectangle in $\sR_i(x)$,
$i\leq m$, can be made as close as we please to the boundary of the corresponding
rectangle of $\sR_i(x_0)$; we are using the continuity of $A(x,h)$ here. It follows that $F(x,u)$ is continuous in $x$,
uniformly over $u$ in compact subsets of $(0,\infty)^2$.

Using Assumption \ref{A.0}(a), the construction also tells us that there
exists $\beta$ such that
\bee\label{Ecompare}
\beta|u|\leq |F(x,u)|\leq \beta^{-1}|u|, \qq x\in \R^d, \quad u\ne 0.
\eee

For each $x$, let $G(x, \cdot)$ be the inverse of $F(x,\cdot)$.
Define
$$N_t(C)=\sum_{s\leq t} 1_{(G(X_{s-}, \Delta X_s)\in C)}$$
and $$\lam(C)=\int_C \frac{du}{|u|^{d+\al}}.$$

\begin{proposition}\label{PA2}
Let $x_0\in \R^d$ and let $\P$ be a solution to the martingale problem
for $\sL$ started at $x_0$.
Then with respect to $\P$, $N_t(\cdot)$ is a Poisson point process with intensity measure $\lam$.
\end{proposition}

\proof Let $F(x,C)=\{F(x,z): z\in C\}$. Then
$$N_t(C)=\sum_{s\leq t} 1_{(\Delta X_s\in F(X_{s-},C))}.$$
By Proposition \ref{P3.21} and a limit argument, the right hand side is equal to a 
martingale plus
$$\int_0^t \int_{F(X_{s-},C)} \frac{A(X_s,h)}{|h|^{d+\al}}\, dh\, ds.$$
By \eqref{Eequiv} and the fact that $X$ has only countably many jumps, this in turn is equal to
\begin{align*}
\int_0^t \int 1_{(F(X_{s-},u)\in F(X_{s-},C))}\frac{du}{|u|^{d+\al}}\, ds
&=\int_0^t \int_C \frac{du}{|u|^{d+\al}}\,ds\\
&=\lam(C)t.
\end{align*}
Therefore by Proposition \ref{PA1} we see that $N_t(\cdot)$ is a Poisson point process.
\qed

Set
$\mu([0,t]\times C)=N_t(C)$.
Note the definition of $\mu$ does not depend on $\P$.

\begin{proposition}\label{AP3}
$X_t$ solves the stochastic differential equation
\begin{align}
X_t&=X_0+\int_0^t \int_{|F(X_{s-},z)|\leq 1} F(X_{s-},z)\, (\mu(dz\, ds)-\lam(dz)\, ds)\nn\\
&~~~ +\int_0^t \int_{|F(X_{s-},z)|> 1} F(X_{s-},z)\, \mu(dz\, ds).\label{EXdef}
\end{align}
\end{proposition}

\proof
Let $\delta>0$. 
Set
$G(x,D)=\{G(x,w): w\in D\}$,  $D^\delta=\{y: |y|>\delta\}$, 
\begin{align*}
H^\delta_t&=\sum_{s\leq t} \Delta X_s 1_{(1\geq |\Delta X_s|>\delta)},\\
\wt{H^\delta}_t
&=\int_0^t \int_{G(X_{s-}, D^\delta\setminus D^1)} F(X_{s-},z)\, \lam(dz)\, ds,\\
K_t&=\sum_{s\leq t} \Delta X_s 1_{( |\Delta X_s|>1)}.
\end{align*}
If $\Delta X_s\ne 0$, the definition of $\mu$ via
$N_t$ shows that $\mu$ assigns unit mass  to some point  $(z,s)$ satisfying 
$z=G(X_{s-},\Delta X_s)$, or with $\Delta X_s=F(X_{s-},z)$. 
Hence
\bee\label{Ag1}
H^\delta_t=
\int_0^t \int_{G(X_{s-}, D^\delta\setminus D^1)}F(X_{s-},z)\, \mu(dz\, ds)
\eee
and
\bee\label{Ag2}
K_t=\int_0^t \int_{G(X_{s-},  D^1)}F(X_{s-},z)\, \mu(dz\, ds).
\eee

By \cite{Jacod} or \cite[Theorem II.10]{LM}, there exists a function 
$\ol F(x,z)$ 
satisfying
\bee\label{Eeq2}
\int 1_B(\ol F(x,u))\frac{du}{|u|^{d+\al}} =\int_B \frac{A(x,h)}{|h|^{d+\al}}
\, dh,
\eee
for $B$ Borel and
and a Poisson point process $\ol \mu$ such that $X_t$ solves
\begin{align*}
X_t&=X_0+\int_0^t \int_{|\ol F(X_{s-},z)|\leq 1} \ol F(X_{s-},z)\, (\ol \mu(dz\, ds)-\lam(dz)\, ds) \\
&~~~+\int_0^t \int_{|\ol F(X_{s-},z)|> 1} \ol F(X_{s-},z)\, \ol \mu(dz\, ds).
\end{align*}
From this equation we see that $\ol \mu$ gives unit mass to a point $(z,s)$
if and only if $\Delta X_s=\ol F(X_{s-},z)$.
It follows that
\begin{align*}
V^\delta_t&=X_t-X_0-K_t-(H^\delta_t-\wt{H^\delta}_t)\\
&=
\int_0^t \int_{|\ol F(X_{s-},z)|\leq \delta} 
\ol F(X_{s-},z) \, (\ol \mu(dz\, ds)-\lam(dz)\, ds).
\end{align*}
A limit argument and \eqref{Eeq2} show that
$$\int |\ol F(x,u)|^2\frac{du}{|u|^{d+\al}}=\int |h|^2\frac{A(x,h)}{|h|^{d+\al}}\,
dh,$$
which is bounded uniformly in $x$. Consequently each component of
$V^\delta$ is a pure jump
martingale and $\E \sup_{s\leq t} |V^\delta_t|^2\to 0$ as $\delta\to 0$.

On the other hand, using \eqref{Ag1} and \eqref{Ag2},

\begin{align*}
V_t^\delta&=
X_t-\Big(X_0+\int_0^t \int_{\delta<|F(X_{s-},z)|\leq 1} 
F(X_{s-},z)\, (\mu(dz\, ds)-\lam(dz)\, ds)\\
&~~~ +\int_0^t \int_{|F(X_{s-},z)|> 1} 
F(X_{s-},z)\, \mu(dz\, ds)\Big).
\end{align*}
Our conclusion follows.
\qed

Define $Y^{n}_s$ to be equal to $x_0$ if $s<1/n$ and equal to
$X_{(k-1)/n}$ if $k/n\leq s<(k+1)/n$.
The reason for the $1/n$ delay will appear in \eqref{E4.11A}.
Let
\begin{align}
X_t^n&=X_0+\int_0^t \int_{|F(Y_{s}^n,z)|\leq 1} F(Y^{n}_s,z)
(\mu(dz\, ds)-\lam( dz)\, ds)\nn\\
&~~~+ \int_0^t \int_{|F(Y_s^n,z)|>1} F(Y^{n}_s,z)\, \mu(dz\, ds).
\label{EXndef}
\end{align}

\begin{proposition}\label{PAlim}
Let $x_0\in \R^d$ and let $\P$ be a solution to the martingale problem
for $\sL$ started at $x_0$.
For each $t_0$
$$\sup_{t\leq t_0} |X_t-X^n_t|\to 0$$
in probability as $n\to 0$.
\end{proposition}

\proof Except for $s=0$, notice 
$Y^{n}_s\to X_{s-}$ a.s.\ under $\P$, using the fact that the
paths of $X_t$ have left limits. 
Except for $z$ in the boundary of any of the $2^d$ orthants, 
$F(Y^{n}_s,z)\to F(X_{s-},z)$ a.s. if $s>0$.

Let $X^{n,I}_t$ be the first double integral on the right hand side
of \eqref{EXndef} and $X^{n,II}_t$ the second. Similarly let
$X^I_t$ be the first double integral on the right in \eqref{EXdef}
and $X^{II}_t$ the second. 
Let
\begin{align*}
Z_t^n&=X_0+\int_0^t \int_{|F(X_{s-},z)|\leq 1} F(Y_s^n,z)
(\mu(dz \, ds)-\lam(dz)\, ds)\\
&~~~~~~+\int_0^t \int_{|F(X_{s-},z)|> 1} F(Y_s^n,z)
\mu(dz \, ds)\\
&=X_0+Z_t^{n,I}+Z_t^{n,II}.
\end{align*}

Using Doob's inequality on each component
and basic properties of stoch\-as\-tic integrals with respect to Poisson
point processes (see, e.g., \cite{Jacod}), 
\begin{align*}
\E\sup_{t\leq t_0} |X^I_t-Z^{n,I}_t|^2&
\leq c\E |X^I_{t_0}-Z^{n,I}_{t_0}|^2\\
&=c\E\int_0^{t_0}\int_{|F(X_{s-},z)|\leq 1}
|F(X_{s-},z)-F(Y^{n}_s,z)|^2 \, \lam(dz)\, ds.
\end{align*}
Using \eqref{Ecompare}, the integrand is bounded by
$$c|z|^2 1_{(|z|\leq \beta^{-1})},$$
which is integrable with respect to $\lam(dz)\, ds$. For $s>0$ and all $z$ not on the boundary of any of the
orthants, and hence for almost every $z$ with respect to $\lam$, the
integrand tends to 0, a.s. Therefore by dominated convergence
$$\E \sup_{t\leq t_0} |X_t^I-Z_t^{n,I}|^2\to 0.$$ 

Since $|F(X_{s-},z)|>1$ implies $|z|\geq \beta$ by \eqref{Ecompare},
with probability one,
there are only finitely many points $(z,s)$ with $|z|\geq \beta$
charged by $\mu$ before time $t_0$. Also with probability one, none
of the $z$ values will lie on the boundary of any of the orthants.
It follows then that $\sup_{t\leq t_0} |X_t^{II}-Z_t^{n,II}|\to 0$ as
$n\to \infty$.

Notice
$$X_t^n -Z_t^n=\int_0^t \int_{C(n,s)} F(Y_s^n,z)\, \lam(dz)\, ds,$$
where
\begin{align*}
C(n,s)&=\{|F(Y_s^n,z)|\leq 1, |F(X_{s-},z)>1\}\\
& ~~~~~~\cup
\{|F(X_{s-},z)|\leq 1, |F(Y_s^n,z)|>1\}.
\end{align*}
Therefore
$$\E \sup_{t\leq t} |X^n_t-Z_n^t|\leq \sum_{i=1}^5 
\E \int_0^{t_0}\int_{D^i(n,s)} |F(Y_s^n,z)|\, \lam(dz)\, ds,$$
where $\gamma>0$ will be chosen in a moment and
\begin{align*}
D^1(n,s)&=\{|F(Y_s^n,z)|\leq 1-\gamma, |F(X_{s-},z)|>1\},\\
D^2(n,s)&=\{|F(Y_s^n,z)|\leq 1, |F(X_{s-},z)|\geq 1+\gamma\},\\
D^3(n,s)&=\{|F(Y_s^n,z)| \geq 1+\gamma, |F(X_{s-},z)|\leq 1\},\\
D^4(n,s)&=\{|F(Y_s^n,z)|> 1, |F(X_{s-},z)|\leq 1-\gamma\},\\
D^5(n,s)&=\{1-\gamma\leq |F(Y_s^n,z)|, |F(X_{s-},z)|<1+\gamma\}.
\end{align*}
By \eqref{Ecompare}, if $|F(X_{s-},z)|\geq 1-\gamma$, then $|z|\geq c$.
So using \eqref{Eequiv} and Assumption \ref{A.0}
\begin{align*}
\E\int_0^{t_0} &\int _{D^5(n,s)} |F(Y_s^n,z)|\, \lam(dz)\, ds\\
&\leq (1+\gamma) \E\int_0^t \int 1_{B(0,1+\gamma)\setminus B(0,1-\gamma)}
(F(X_{s-},z))\, \lam(dz)\, ds\\
&=(1+\gamma) \E\int_0^{t_0} \int_{B(0,1+\gamma)\setminus B(0,1-\gamma)}
\frac{A(X_{s-},h)}{|h|^{d+\al}}\, dh\, ds\\
&\leq c\gamma t_0.
\end{align*}
So the integral over $D^5(n,s)$ can be made as small as we like by
taking $\gamma$ sufficiently small. Once $\gamma$ is chosen, observe
that $1_{D^1(n,s)}\to 0$ a.s. for every $s>0$ because $Y^n_s\to X_{s-}$. 
Also, on $D^1(n,s)$, we have $|F(X_{s-},z)|>1$, and as above $|z|>c$,
so $|F(Y_s^n,z)|1_{D^1(n,s)}$ is dominated by $(1+\gamma) 1_{(|z|\geq c)}$,
which is integrable with respect to $\lam(dz)\, ds$. So by 
dominated convergence,
$$\E \int_0^{t_0} \int _{D^1(n,s)}|F(Y_s^n,z)|\, \lam(dz)\, ds\to 0.$$
The argument for $D^2(n,s), D^3(n,s)$, and $D^4(n,s)$ is the same.
Hence
$$\E\sup_{t\leq t_0} |X^n_t-Z^n_t|\to 0$$.
\qed

\begin{proposition}\label{Papprox}
Let $x_0\in \R^d$ and let $\P$ be a solution to the martingale problem
for $\sL$ started at $x_0$.
If $f\in C^2_b$,
then
$$f(X^n_t)-f(X^n_0)-\int_0^t \sM^{Y^{n}_s} f(X^n_{s})\, ds$$
is a martingale under $\P$,
where $\sM^y$ is defined in \eqref{DL0}.
\end{proposition}

\proof
If $\mu$ assigns unit mass to $(z,s)$, then $\Delta X_s^n=F(Y^{n}_s,z)$.
By Ito's formula
\begin{align*}
f(&X_t^n)-f(X_0^n)\\
&=\mbox{ martingale }
+\int_0^t \int_{|F(X_{s-},z)|>1} \grad f(X_s^n)\cdot F(Y^{n}_s,z)\, \mu(dz\, ds)\\
&~~~~~~~+\sum_{s\leq t}[f(X_s^n)-f(X^n_{s-})-\grad f(X_{s-}^n)\cdot \Delta X^n_s]\\
&=\mbox{ martingale }+\int_0^t \int [f(X^n_{s-}+F(Y^{n}_s,z))-f(X^n_{s-})\\
&~~~~~~~~~-\grad f(X^n_{s-})\cdot F(Y^{n}_s,z)
1_{(|F(Y^{n}_s,z)|\leq 1)}] \, \mu(dz\, ds)\\
&=\mbox{ martingale }+\int_0^t \int [f(X^n_{s-}+F(Y^{n}_s,z))-f(X^n_{s-})\\
&~~~~~~~~~-\grad f(X^n_{s-})\cdot F(Y^{n}_s,z)
1_{(|F(Y^{n}_s,z)|\leq 1)}] \, |z|^{-(d+\al)}\, dz\, ds.
\end{align*}
Fix $y$ and if $\al\geq 1$,  let $$g(v)=f(y+v)-f(y)-\grad f(y)\cdot v 1_{(|v|\leq 1)}.$$
A limit argument using \eqref{Eequiv} shows
$$\int g(F(x,z))\frac{1}{|z|^{d+\al}}\, dz=\int g(h)\frac{A(x,h)}{|h|^{d+\al}}\, dh.$$
Now taking $y=X^n_{s-}$ and $x=Y^{n}_s$ shows that
$f(X^n_t)-f(X^n_0)$ is equal to a martingale plus
$$\int_0^t \int [f(X_{s-}^n+h)-f(X_{s-}^n)
-\grad f(X_{s-}^n)\cdot h 1_{(|h|\leq 1)}]
\frac{A(Y^{n}_s,h)}{|h|^{d+\al}}\, dh\, ds,$$
which proves the proposition when
$\al\geq 1$. The case $\al<1$ is similar.
\qed

\section{Existence and uniqueness}\label{sect:exist}

\begin{theorem}\label{T4.2}
Suppose Assumption \ref{A.0} holds. Then for each $x$
there exists a solution to the martingale problem for $\sL$ started at
$x$.
\end{theorem}

\proof In view of Propositions \ref{P3.1A} and \ref{P3.3}, existence
of a solution follows by the proof in \cite[Section 3]{Upjump}, with 
minor modifications to handle the case of $d$ dimensions.
\qed

\begin{remark}\label{rem:exist}
{\rm It is easy to see by the same arguments that existence holds if $A(x,h)$
is bounded above and below by positive constants and 
for each $h$,  $A(x,h)$ is continuous in $x$.
}
\end{remark}

We now turn to the proof of uniqueness.
Fix $x_0\in \R^d$. If $\sG$ is the set of solutions to the
martingale problem for $\sL$ started at $x_0$, then $\sG$ is a tight family by the proof
in \cite[Section 3]{Upjump}. Any subsequential limit point of $\sG$ is in
$\sG$ by the arguments in that  same section, and therefore $\sG$ is compact.
Hence by the proofs in \cite[Chapter 12]{StrVar}, it suffices to consider
uniqueness of
strong Markov families of solutions $\{\P^x\}$ to
the martingale problem for $\sL$.

We will sometimes make the following temporary assumption, where we will choose $\zeta$ later:

\ms
\begin{assumption}\label{A.1}
 There exists $\zeta$  such that  
$$|A(x,h)-A(x_0,h)|\leq \frac{\zeta}{\psi_\eta(|h|)}, \qq x\in \R^d, \quad |h|\leq 1.$$

\end{assumption}
\ms

For the rest of this section we take $\al\geq 1$, the case $\al<1$ being
similar.

Let $\sM^z f$ be defined by \eqref{DL0}
and
let $R^z_\lam$ be the corresponding resolvent.
Define an operator $\sH$ by
\begin{align*}
\sH f(x)&= \int_{|h|\leq 1} |f(x+h)-f(x)-\grad f(x)\cdot h|\frac{\zeta}{\psi_\eta (|h|) |h|^{d+\al}}\, dh\\
&~~~~+\int_{|h|>1} |f(x+h)-f(x)|\, \frac{dh}{|h|^{d+\al}}, \qq f\in C^2_b.
\end{align*}

\begin{proposition}\label{PEU1}
There exists a constant $c_1$ not depending on $x_0$  such that
\bee\label{E200}
\norm{\sH R^{x_0}_\lam f}_2\leq c_1(\zeta+\lam^{-1})\norm{f}_2, \qq f\in L^2\cap C^2_b.
\eee
\end{proposition}

\proof By Minkowski's inequality for integrals, 
\begin{align}
\|\sH R^{x_0}_\lam f\|_2&\leq
\int_{|h|\leq 1} \norm{R^{x_0}_\lam f(x+h)-R^{x_0}_\lam f(x) \nn\\
&\qq\qq\qq-\grad R^{x_0}_\lam f(x)\cdot h}_2 \frac{\zeta}{\psi_\eta(h)|h|^{d+\al}}\, dh\nn\\
&~~~+
\int_{|h|>1} \norm{R^{x_0}_\lam f(x+h)-R^{x_0}_\lam f(x) }_2 \frac{c}{|h|^{d+\al}}\, dh.
\label{E21}
\end{align}
By Proposition \ref{P3.8}(c), the first term on the right of
\eqref{E21} is bounded by
\bee\label{Term1}
c\int_{|h|\leq 1} |h|^\al \frac{\zeta}{\psi_\eta(h)|h|^{d+\al}}\, dh
\norm{f}_2\leq c\zeta \norm{f}_2.\eee
By Proposition \ref{P3.8}(b) the second term on the right of
\eqref{E21} is bounded by
\bee\label{Term3}
\frac{c}{\lam}\int_{|h|>1} \frac{dh}{|h|^{d+\al}}\, \norm{f}_2.\eee
\qed

\begin{corollary}\label{VC1}
Suppose Assumption \ref{A.1} holds.
There exists $\kappa$ such that
\bee\label{VE1}
\norm{(\sL-\sM^{x_0})R_\lam^{x_0}f}_2\leq \kappa(\zeta+\lam^{-1})\norm{f}_2, 
\qq f\in L^2\cap C^2_b,
\eee
and
\bee\label{VE2}
\norm{\sup_{w\in \R^d} |\sM^w R_\lam^{x_0} f(\cdot)-\sM^{x_0}R_\lam^{x_0}
f(\cdot)}_2\leq \kappa(\zeta+\lam^{-1})\norm{f}_2, \qq 
f\in L^2\cap C^2_b.
\eee
\end{corollary}

\proof
If Assumption \ref{A.1} holds, then
\begin{align}
|(\sL -\sM^{x_0}) &R^{x_0}_\lam f(x)|\nn\\
&=\Big|\int [R^{x_0}_\lam f(x+h)-R^{x_0}_\lam f(x)\label{E201}\\
&~~~~~~~-\grad R^{x_0}_\lam f(x)\cdot h 1_{|h|\leq 1)}]
\frac{A(x,h)-A(x_0,h)}{|h|^{d+\al}}\, dh\Big|\nn\\
&\leq c\sH R^{x_0}_\lam f(x)\nn
\end{align}
and for each $w$
\begin{align}
|(\sM^w R^{x_0}_\lam &f(x,y) -\sM^{x_0} R^{x_0}_\lam f(x)|\nn\\
&\leq \Big|\int_{|h|\leq 1} [R^{x_0}_\lam f(x+h)-R^{x_0}_\lam f(x)\label{E202}\\
&~~~~~~~-\grad R^{x_0}_\lam f(x)\cdot h 1_{|h|\leq 1)}]
\frac{A(w,h)-A(x_0,h)}{|h|^{d+\al}}\, dh\Big|\nn\\
&~~~~~~~+\Big|\int_{|h|\geq 1}
[R_\lam^{x_0} f(x+h)-R^{x_0}_\lam f(x)] \frac{A(w,h)+A(x_0,h)}
{|h|^{d+\al}}\, dh\Big|\nn\\
&\leq c\sH R^{x_0}_\lam f(x).\nn
\end{align}
Now combine Proposition \ref{PEU1}, \eqref{E201}, and \eqref{E202}.
\qed

\begin{proposition}\label{VP1}
Let $\{\P^x\}$ be a strong Markov family of solutions to the martingale
problem for $\sL$. Set
$$S_\lam f(x)=\E^x \int_0^\infty e^{-\lam t} f(X_t)\, dt.$$
Suppose Assumption \ref{A.1} holds with $\zeta$ and $\lam$ chosen
so that $\kappa(\zeta+\lam^{-1})\leq 1/2$, where $\kappa$ is as
in Corollary \ref{VC1}. Let $\rho\in L^2$ be non-negative
with compact support. Then 
$$\sup_{\norm{g}_2\leq 1} \Big|\int S_\lam g(x)\rho(x)\, dx\Big|<\infty.$$
\end{proposition}

\proof Define $X_t^n$ as in Section \ref{sect:approx} and define
$$S_\lam^n g(x)=\E^x\int_0^\infty e^{-\lam t} g(X_t^n)\, dt, \qq g\in C^2_b.$$

\ni{\sl Step 1:} Our first goal is to show 
that if
\bee\label{VE3}
\Lambda_n=\sup_{\norm{g}_2\leq 1}\Big|\int S_\lam^n g(x)\rho(x)\, dx\Big|,
\eee
then
$\Lambda_n<\infty$.
The value of $\Lambda_n$ will depend on $\rho$.

To prove \eqref{VE3} it suffices to suppose $g\geq 0$ since we can write 
an arbitrary $g$ as the difference of its positive and negative parts.
Suppose $g\in C^2_K$ and write
\bee\label{VE4}
S_\lam^ng(x)=\E^x\int_0^{1/n}e^{-\lam t} g(X_t^n)\, dt
+\sum_{k=1}^\infty \E^x \int_{k/n}^{(k+1)/n} e^{-\lam t} g(X_t^n)\, dt.
\eee
Over the time interval $[0,1/n)$, the process $X^n_t$ behaves like the L\'evy process
corresponding to $\sM^{x_0}$ started at $x$. So the first term on the
right hand side of \eqref{VE4} is bounded by
$R^{x_0}_\lam g(x)$. By the Cauchy-Schwarz inequality
and Proposition \ref{P3.5},
\bee\label{VE5}\Big|\int R_\lam ^{x_0} g(x)\rho(x)\, dx\Big|
\leq \norm{R_\lam^{x_0}g}_2\norm{\rho}_2\leq \frac{c}{\lam} \norm{g}_2.
\eee

The $k^{th}$ term on the right hand side of \eqref{VE4} is 
\begin{align*}
e^{-\lam(k-1)/n}&\E^x\int_{k/n}^{(k+1)/n} e^{-\lam (t-(k-1)/n)} g(X_t^n)\, dt\\
&\leq ce^{-\lam k/n}\E^x\Big[\E^x\Big[\int_{1/n}^{2/n} g(X^n_{t+\frac{k-1}{n}})
\, dt\mid \F_{(k-1)/n}\Big]\, \Big].
\end{align*}
Let us temporarily write $\ol Y$ for $Y^n_{(k-1)/n}$. Conditional on
$\F_{(k-1)/n}$, the process $X_t^n$ over the time interval $[k/n, (k+1)/n)$
behaves like the L\'evy process corresponding to $\sM^{\ol Y}$ started
at $X^n_{(k-1)/n}$ and run over the time interval $[1/n,2/n]$. Therefore
\begin{align}
\E^x\Big[\int_{1/n}^{2/n} &g(X^n_{t+\frac{k-1}{n}})
\, dt\mid \F_{(k-1)/n}\Big]\label{E4.11A}\\
&\leq \int_{1/n}^{2/n} P_t^{\ol Y} g(X^n_{(k-1)/n})\, dt\nn\\
&\leq e^{\lam/n}P^{\ol Y}_{1/n} R_\lam^{\ol Y} g(X^n_{(k-1)/n}).\nn
\end{align}
Using Corollary \ref{C3.7} we have
\begin{align*}
\P^w_{1/n} R_\lam^w g(v)&=\int \ol p^w(1/n,v-z)R_\lam^wg(z)\, dz\\
&\leq \norm{\ol p^w(1/n, \cdot)}_2\norm{R^w_\lam g}_2\\
&\leq c_n \frac{1}{\lam} \norm{g}_2,
\end{align*}
where $c_n$ depends on $n$.
Hence the $k^{th}$ term on the right hand side of \eqref{VE4}
is bounded by $ce^{-\lam k/n}\norm{g}_2$. Because $\rho$ is in $L^2$ with
compact support, 
$$\int \E^x\int_{k/n}^{(k+1)/n} g(X_t^n)\, dt\, \rho(x)\, dx
\leq c\norm{g}_2\int \rho(x)\, dx\leq c\norm{g}_2.$$
Combining this with \eqref{VE5} and taking the supremum over
$g\in C^2_K$ with $\norm{g}_2\leq 1$  proves \eqref{VE3}.
\ms

\ni{\sl Step 2:} Next we show there exists a constant $\Lambda<\infty$
independent of $n$ such that
\bee\label{VE6}
\sup_{\norm{g}_2\leq 1} \Big|\int S^n_\lam g(x) \rho(x)\, dx\Big|\leq \Lambda.
\eee

Let $f\in C^2_b$. By Proposition \ref{Papprox}
$$\E^x f(X_t^n)-f(x)=\E^x\int_0^t \sM^{Y^n_s}f(X^n_s)\, ds.$$
Multiplying by $e^{-\lam t}$ and integrating over $t$ from 0 to $\infty$
\begin{align}
S_\lam^nf(x)-\frac{1}{\lam}f(x)&=\E^x\int_0^\infty e^{-\lam t}
\int_0^t \sM^{Y^n_s}f(X_s^n)\, ds\, dt\label{VE7}\\
&=\E^x\int_0^\infty \sM^{Y^n_s}f(X^n_s)\int_s^\infty e^{-\lam t}\, dt\, ds\nn\\
&=\frac{1}{\lam}\E^x\int_0^\infty e^{-\lam s} \sM^{Y^n_s}f(X^n_s)\, ds\nn\\
&=\frac{1}{\lam} S_\lam^n \sM^{x_0}f(x)+\frac{1}{\lam}\E^x\int_0^\infty e^{-\lam s}(\sM^{Y^n_s}-\sM^{x_0})f(X_s^n)\, ds.\nn
\end{align}
If $g\in C^2_K$, set $f=R_\lam^{x_0}g$. By translation invariance, $f\in C^2_b$.
Standard semigroup manipulations show 
$$\sM^{x_0} f=\sM^{x_0}R^{x_0}_\lam g=\lam R^{x_0}_\lam g-g.$$
Therefore
$$S^n_\lam R^{x_0}_\lam g(x)-\frac{1}{\lam} R^{x_0}_\lam g(x)
\leq S^n_\lam R_\lam^{x_0} g(x)-\frac{1}{\lam}S^n_\lam g(x)+
\frac{1}{\lam} S_\lam^n H(x),$$
where
$$H(y)=\sup_{w\in \R^d}|\sM^wR^{x_0}_\lam g(y)-\sM^{x_0}R^{x_0}_\lam g(x)|.$$
We thus have
\bee\label{VE8}
S^n_\lam g(x)\leq R^{x_0}_\lam g(x)+S^n_\lam H(x).
\eee
By Corollary \ref{VC1}
and our  choice of $\zeta$ and $\lam$,
$$\norm{H}_2\leq \kappa(\zeta+\lam^{-1})\norm{g}_2\leq \frac12 \norm{g}_2.$$

Multiplying \eqref{VE8} by $\rho(x)$ and integrating, 
\begin{align*}
\Big|\int S^n_\lam g(x)\rho(x)\, dx\Big|&\leq
\Big|\int R^{x_0}_\lam g(x)\rho(x)\, dx\Big|
+\Big|\int S^n_\lam H(x)\rho(x)\, dx\Big|\\ 
&\leq \norm{R_\lam^{x_0}g}_2\norm{\rho}_2+\Lambda_n \norm{H}_2\\
&\leq \frac{1}{\lam}\norm{\rho}_2\norm{g}_2+\frac12\Lambda_n \norm{g}_2,
\end{align*}
where $\Lambda_n$ is defined in Step 1.
Taking the supremum over $g\in C^2_K$ with $\norm{g}_2\leq 1$, we thus
have
$$ \Lambda_n\leq \frac{\norm{\rho}_2}{\lam}+\frac12 \Lambda_n.$$
In Step 1 we proved $\Lambda_n<\infty$, and we conclude
$$\Lambda_n\leq \frac{2}{\lam}\norm{\rho}_2.$$
\ms

\ni{\sl Step 3:} We now pass to the limit in $n$.
By Step 1 and Step 2, if $g\in C^2_K$ with $\norm{g}_2\leq 1$,
then 
$$\Big|\int S^n_\lam g(x)\rho(x)\, dx\Big|\leq \frac{2\norm{\rho}_2}{\lam}.$$
On the other hand,
$$S_\lam^n g(x)=\E^x \int_0^\infty e^{-\lam t} g(X_t^n)\, dt
\to \E^x \int_0^\infty e^{-\lam t}g(X_t)\, dt=S_\lam g(x)$$
by dominated convergence. We thus see that
$$\Big|\int S_\lam g(x)\rho(x)\, dx\Big|\leq \frac{2\norm{\rho}_2}{\lam}.$$
Our result follows by taking the supremum over $g\in C^2_K$ with $\norm{g}_2\leq 1$.
\qed

\ni{\bf Proof of Theorem \ref{main}:} Let $x_0\in \R^d$. Let $\rho\in L^2$
with compact support. We have seen that it suffices to prove uniqueness
when we have a strong Markov family of solutions to the martingale
problem for $\sL$, so suppose we have two such families $\{\P^x_i\}$, $i=1,2$.
Define 
$$S_\lam^i f(x)=\E^x_i \int_0^\infty e^{-\lam t}\, dt, \qq i=1,2,$$
and let
$$S^\Delta_\lam =S_\lam^1-S_\lam^2.$$
Suppose $\lam_0$ and $\zeta$ are chosen so that $\kappa(\zeta+\lam_0^{-1})
\leq \frac12$ and $\lam>\lam_0$, and suppose Assumption \ref{A.1} holds
with this choice of $\zeta$.

Since $\P^x_i$ is a solution to the martingale problem for $\sL$ started at
$x$, for $f\in C^2_b$
$$\E^x_i f(X_t)-f(x)=\E^x_i \int_0^t \sL f(X_s)\, ds.$$
Multiplying by $e^{-\lam t}$ and integrating over $t$ from 0 to $\infty$,
\begin{align*} S_\lam^if(x)-\frac{1}{\lam}f(x)
&=\E^x_i \int_0^\infty e^{-\lam t} \int_0^t \sL f(X_s)\, ds\, dt\\
&=\E^x_i \int_0^\infty \sL f(X_s)\int_s^\infty e^{-\lam t}\, dt\, ds\\
&=\frac{1}{\lam} S^i_\lam \sL f(x)\\
&=\frac{1}{\lam}S^i_\lam \sM^{x_0}f(x)+ 
\frac{1}{\lam}S^i_\lam (\sL-\sM^{x_0})f(x).
\end{align*}
Now take $g\in C^2_K$ and set $f=R^{x_0}_\lam g$. Then $f\in C^2_b$ and $\sM^{x_0}f
=\lam R^{x_0}_\lam g-g$.
Hence
$$S_\lam^i R_\lam^{x_0}g(x)
-\frac{1}{\lam} R^{x_0}_\lam g(x)=S^i_\lam R^{x_0}_\lam g(x)-\frac{1}{\lam}S_\lam^g(x)
+\frac{1}{\lam} S^i_\lam (\sL-\sM^{x_0})R_\lam^{x_0}g(x),$$
or
\bee\label{VE9}
S_\lam^ig(x)=R_\lam^{x_0} g(x)+S_\lam^i(\sL-\sM^{x_0})R^{x_0}_\lam g(x).
\eee

Let $$\Theta=\sup_{\norm{g}_2\leq 1}\Big| \int S_\lam^\Delta g(x)\rho(x)\, dx\Big|.$$
By Proposition \ref{VP1}, we know that $\Theta<\infty$. From \eqref{VE9}
$$S_\lam^\Delta g(x)=S_\lam^\Delta (\sL-\sM^{x_0})R_\lam^{x_0}g(x).$$
Multiplying by $\rho(x)$ and integrating,
\begin{align*}
\Big|\int S_\lam^\Delta g(x)\rho(x)\, dx\Big|
&=\Big|\int S_\lam^\Delta (\sL-\sM^{x_0})R_\lam^{x_0}g(x)\rho(x)\, dx\Big|\\
&\leq \Theta \norm{(\sL-\sM^{x_0})R_\lam ^{x_0} g}_2.
\end{align*}
By Corollary \ref{VC1} this is bounded by $\frac12 \norm{g}_2$. Taking the supremum 
over $g\in C^2_K$ with $\norm{g}_2\leq 1$, we then obtain $\Theta\leq \frac12
\Theta$. Since $\Theta<\infty$, this implies $\Theta=0$.
This can be rewritten as
$$\int S_\lam^1g(x)\rho(x)\, dx=\int S_\lam^2 g(x) \rho(x)\, dx.$$

This is true for each  $\rho \in L^2$ with compact support, and we conclude $S_\lam^1g(x)=
S_\lam^2 g(x)$ for almost every $x$. By Proposition \ref{P3.4},
$S_\lam^ig(x)$ is continuous in $x$, so we have equality for all $x$.
By the uniqueness of the Laplace transform and the right continuity of
$X_t$, we conclude
$$\E^x_1 g(X_t)=\E^x_2 g(X_t)$$
for all $x$ and all $t$ whenever $g$ is continuous and bounded. By a limit argument this equality holds
for all bounded $g$. Finally, by using the Markov property, the
finite dimensional distributions under $\P^x_1$ and $\P^x_2$ are
the same for each $x$.

The last step
is to remove the use of Assumption \ref{A.1}. This is a standard
localization argument. Because of Assumption \ref{A.0},
there exists $\wt A(x,h)$ such that $\wt A$ agrees with $A$ in a 
neighborhood of $x_0$ and such that Assumption \ref{A.0}
holds for $\wt A$. If $\wt \sL$ is the operator defined in terms
of $\wt A$ in the same way as $\sL$ is defined in terms of $A$, the
above shows we have uniqueness for the martingale problem for
$\wt \sL$ started at $x_0$. From this point on, we proceed exactly
as in the diffusion case;  see \cite[Chapter VI]{DEO}.
This completes the proof of Theorem \ref{main}.
\qed

\begin{remark}\label{rem:psi}
{\rm
It is clear that $\psi_\eta$ can be replaced by any decreasing function $\psi$ such
that 
$$\int_{|h|\leq 1} \frac{1}{\psi(|h|)|h|^{d+\al}}\, dh<\infty.$$
}
\end{remark}

\begin{remark}\label{rem:lack}
{\rm 
Just as in the case of diffusions, we do not really need continuity
of $A(x,h)$ in $x$, just that each point $x_0$ has a neighborhood
in which $\ol A(x, \cdot)$ is sufficiently close to $\ol A(x_0,\cdot)$.
}
\end{remark}

\begin{remark}\label{rem:kom}
{\rm
In \cite{komatsu} Komatsu considers uniqueness for operators of the form
$\sL_1+\sL_2$, where $\sL_1$ is a stable process of index $\al$
(not necessarily symmetric, but he requires that the jump kernel
for $\sL_1$ be $d$ times continuously differentiable in $h$ away
from the origin) and 
$$\sL_2 f(x)=\int [f(x+h)-f(x)] n(x,dh)$$
(with the appropriate modification when $\al\geq 1$) 
where there exists a measure $n^*$ such that $|n(x,dh)|\leq n^*(dh)$
and $\int (1\land |h|^\al) \, n^*(dh)<\infty$. If we write
the kernel for $\sL_1$ as $A_0(h)/|h|^{d+\al}$ and if in
addition we assume $n^*$ has a density
with respect to Lebesgue measure, we can fit his framework into ours
by setting $$A(x,h)=A_0(h)+\frac{n(x,dh)}{dh}|h|^{d+\al}.$$
}
\end{remark}

\begin{remark}\label{rem:largejumps}
{\rm 
We have not tried to find the weakest possible conditions possible,
particularly with regard to ``large jumps.'' One will still have
uniqueness with minimal assumptions on the intensity of the jumps
above some size $\delta$. This is apparent from the stochastic differential
equations representation of $X$: there are only finitely many jumps
of size larger than $\delta$ in any finite time interval, and so one can
consider them sequentially. Our results will imply uniqueness up to
the time of the first jump of size larger than $\delta$, the law of
that jump is uniquely determined by the location the process jumps from,
and one then has uniqueness up to the time of the second large jump, and so on. 
}
\end{remark}

%ZZZZZZZZZZZZZZZZZZZZZZZZZZZZZZZZ

\medskip

\begin{minipage}[t]{0.39\textwidth}
{\bf Richard F. Bass}\\
Department of Mathematics\\
University of
Connecticut \\
Storrs, CT 06269-3009, USA\\
{\it bass@math.uconn.edu}
\end{minipage}
\hfill
\begin{minipage}[t]{0.55\textwidth}
{\bf Huili Tang}\\
Department of Mathematics\\
University of Connecticut\\
Storrs, CT 06269-3009, USA\\
{\it huili@math.uconn.edu} 
\end{minipage}

\end{document}